\documentclass[11pt]{article}

\usepackage{tikz}
\usepackage{subfigure}
\usepackage[english]{babel}
\usepackage{graphicx}

\usepackage[center]{caption2}
\usepackage{amsfonts,amssymb,amsmath,latexsym,amsthm}
\usepackage{multirow}

\def\ex{\mbox{ex}}
\def\sat{\mbox{sat}}
\def\Sat{\mbox{Sat}}

\newtheorem{thm}{Theorem}[section]
\newtheorem{cor}[thm]{Corollary}
\newtheorem{lem}[thm]{Lemma}
\newtheorem{obs}[thm]{Observation}

\newtheorem{conj}[thm]{Conjecture}

\newtheorem{claim}{Claim}
\newtheorem{fact}{Fact}

\begin{document}

\title{Minimum saturated graphs without $4$-cycles and $5$-cycles
\thanks{The work was supported by National Natural Science Foundation of China (No.12401455) to Yue Ma.}
}
\author{Yue Ma$^a$\\
\small $^{a}$School of Mathematics and Statistics,\\
\small Nanjing University of Science and Technology,\\
\small Nanjing, Jiangsu 210094, China.\\
\small $^a$yma@njust.edu.cn
}

\date{}

\maketitle

\begin{abstract}
Given a family of graphs $\mathcal{F}$, a graph $G$ is said to be $\mathcal{F}$-saturated if $G$ does not contain a copy of $F$ as a subgraph for any $F\in\mathcal{F}$, but the addition of any edge $e\notin E(G)$ creates at least one copy of some $F\in\mathcal{F}$ within $G$.
The minimum size of an $\mathcal{F}$-saturated graph on $n$ vertices is called the saturation number, denoted by $\sat(n, \mathcal{F})$. 
Let $C_r$ be the cycle of length $r$.
In this paper, we study on $\sat(n, \mathcal{F})$ when $\mathcal{F}$ is a family of cycles. In particular, we determine that $\sat(n, \{C_4,C_5\})=\lceil\frac{5n}{4}-\frac{3}{2}\rceil$ for any positive integer $n$.
\end{abstract}

\section{Introduction}
Let $G=(V,E)$ be a graph. We call $|V|$ the {\it order} of $G$ and $|E|$ the {\it size} of it. If $|V|=n$, we call $G$ an $n$-vertex graph.
For positive integer $r$ and $s$, let $K_r$ denote the complete graph on $r$ veritces and $K_{r,s}$ denote the complete bipartite graph whose two partition sets consist of $r$ vertices and $s$ vertices respectively.
For any $v\in V(G)$, we use $N_G(v)$ to denote the set of the neighbors of $v$ in $G$ and let $d_{G}(v)=|N_{G}(v)|$ be the {\it degree} of $v$ in $G$. 
For a finite set $U$, let $K[U]$ be the complete graph on $U$. For any two sets $U_1$ and $U_2$, let $K[U_1,U_2]$ be a complete bipartite graph with two partition sets $U_1$ and $U_2$.
For a set $U\subseteq V(G)$, let $G[U]=G\cap K[U]$ be the subgraph induced by $U$. 
Write $G-U=G[V(G)\backslash U]$. 
Write $N_{G}(U)=\bigcup_{v\in U}N_{G}(v)$.
For any two disjoint sets $U_1,U_2\subset V(G)$, let $G[U_1,U_2]=G\cap K[U_1,U_2]$.
For integers $a$ and $b$ with $a\le b$, let $[a,b]=\{a,a+1,a+2,\dots, b\}$ and $[a,+\infty)=\{r\in\mathbb{Z}:r\ge a\}$. 
Let$P_r$ denote the path of length $r-1$ (on $r$ vertices). Let $C_r$ denote the cycle of length $r$ for $r\ge 3$.
We use $x_1x_2\dots x_r$ to denote the path with edge set $\{x_1x_2,x_2x_3,\dots, x_{r-1}x_r\}$ and use $x_1x_2\dots x_rx_1$ to denote the cycle with edge set $\{x_1x_2, x_2x_3, \dots, x_{r-1}x_r, x_rx_1\}$.
For any subset $I\subseteq[3,+\infty)$, let $\mathcal{C}_{I}$ be the family of all cycles of length $r\in I$. 

Given a family $\mathcal{F}$ of graphs, a graph $G$ is said to be {\em $\mathcal{F}$-saturated} if $G$ does not contain a subgraph isomorphic to any member  $F\in\mathcal{F}$ but  $G+e$ contains at least one copy of some $F\in\mathcal{F}$ for any edge $e\notin E(G)$.
The famous {\it Tur\'{a}n number} $\ex(n,\mathcal{F})$ of $\mathcal{F}$ is the maximum number of edges in an $n$-vertex $\mathcal{F}$-saturated graph.
In this paper, we consider the minimum number of edges in an $n$-vertex $\mathcal{F}$-saturated graph, which is called the {\it saturation number}, denoted by $\sat(n,\mathcal{F})$, i.e.
$$\sat(n,\mathcal{F})=\min\{|E(G)| : G\mbox{ is an $n$-vertex }\mathcal{F}\mbox{-saturated graph}\}\mbox{.}$$
We call an $n$-vertex $\mathcal{F}$-saturated graph of size $\sat(n,\mathcal{F})$ a {\it minimum extremal graph} for $\mathcal{F}$ and let $\Sat(n,\mathcal{F})$ be the family of all $n$-vertex minimum extremal graphs for $\mathcal{F}$.


For a single cycle $C_r$, there are many known results for $\sat(n, C_r)$, such as:
\begin{itemize}
\item (Erd\H{o}s, Hajnal and Moon~\cite{EHM64})  $\sat(n,C_3)=n-1$ for $n\ge 3$;

\item (Ollmann~\cite{Oll72}, Tuza~\cite{Tuz89}, Fisher et al.~\cite{FFL97})  $\sat(n, C_4)=\lfloor\frac{3n-5}{2}\rfloor$ for $n\ge 5$;

\item (Chen~\cite{Che09,Che11}) $\sat(n, C_5)=\lceil\frac{10}{7}(n-1)\rceil$ for $n\ge 21$;

\item(Lan, Shi, Wang and Zhang~\cite{Lan21}) $\frac{4n}{3}-2\le\sat(n,C_6)\le\frac{4n+1}{3}$ for $n\ge 9$;

\item (F\"{u}redi and Kim~\cite{FK13})  $(1+\frac{1}{r+2})n-1<\sat(n,C_r)<(1+\frac{1}{r-4})n+\binom{r-4}{2}$ for all $r\ge 7$ and $n\ge 2r-5$;

\item
(Clark, Entringer and Shapiro~\cite{Clark83, Clark92}, Lin et al.~\cite{LJZY97})  $\sat(n,C_n)=\lceil\frac{3n}{2}\rceil$ for $n=17$ or $n\ge 19$.
\end{itemize}

For cycle families with at least two cycles, there are not many works been done. Here are some known exact results for $\mathcal{C}_{[r,+\infty)}$.
\begin{itemize}
\item (Ferrara et al.~\cite{Subdivision12}) $\sat(n, \mathcal{C}_{[4,+\infty)})=\lceil\frac{5n}{4}-\frac{3}{2}\rceil$ for $n\ge 1$.

\item (Ferrara et al.~\cite{Subdivision12}) $\sat(n, \mathcal{C}_{[5,+\infty)})=\lceil\frac{10(n-1)}7\rceil$ for $n\ge 5$.

\item (Ma, Hou, Hei and Gao~\cite{Our21}) $\sat(n, \mathcal{C}_{[6,+\infty)})=\lceil\frac{3(n-1)}2\rceil$ for $n\ge 10$.

\end{itemize}

For a survey for other saturated results, we refer to Faudree, Faudree, and Schmitt~\cite{FF11}.
In this paper, we study on the saturation number of $\mathcal{C}_I$ with $|I|\ge 2$. 

One may think it is natural to consider $I=\{3,4\}$ at first. However, the following observation makes this trivial.
\begin{obs}\label{3trivial}
For an integer $n\ge 1$ and an integer set $I$ with $3\in I\subseteq[3,+\infty)$, $\sat(n,\mathcal{C}_{I})=n-1$. 
\end{obs}
\begin{proof}
Note that a $\mathcal{C}_{I}$-saturated graph $G$ on $n$ vertices must be connected. Otherwise, the addition of any non-edge between two components of $G$ would create no new cycles. Thus, $\sat(n,\mathcal{C}_{I})\ge n-1$.
On the other hand, since the addition of any non-edge $e$ of $H=K_{1,n-1}$ would create a $C_3\in \mathcal{C}_I$ in $H+e$, and $H$ itself does not contain any cycles, we can see $H$ is a $\mathcal{C}_{I}$-saturated graph $G$ on $n$ vertices. Thus, $\sat(n,\mathcal{C}_{I})\le n-1$. This completes the proof.
\end{proof}

By Observation~\ref{3trivial}, the first nontrivial case should be $I=\{4,5\}$. The following is the main result of this paper, which gives the exact value of $\sat(n,\mathcal{C}_{I})$ when $I=\{4,5\}$.

\begin{thm}\label{45}
For $n\ge 1$, $\sat(n,\mathcal{C}_{\{4,5\}})=\lceil\frac{5n}{4}-\frac{3}{2}\rceil$.
\end{thm}

The rest of the paper is arranged as follows. We give the proof of Theorem~\ref{45} in Section 2. In Section 3, there are some remarks and conjectures for $\mathcal{C}_{I}$-saturated graphs. 

\section{Proof of Theorem~\ref{45}}
\subsection{Upper bound: a construction for a $\mathcal{C}_{\{4,5\}}$-saturated graph}
\begin{figure}[h]
\centering
\includegraphics[width=5in]{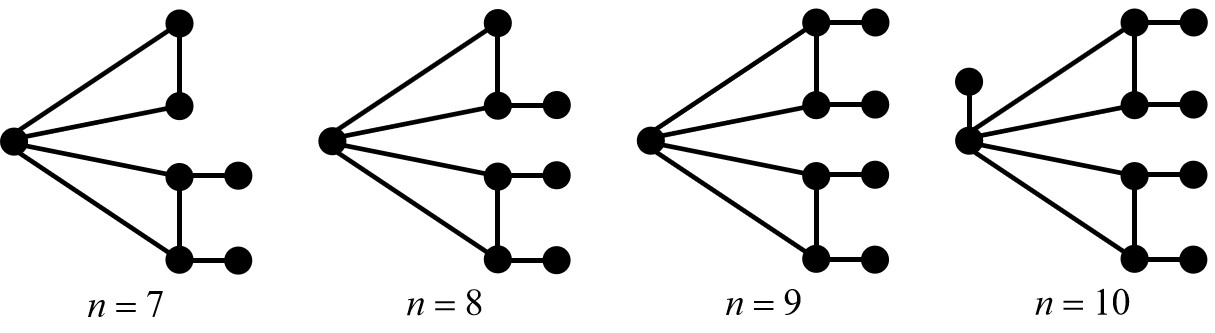}
\caption{$Sat_n$ for $n\in[7,10]$}\label{sat4}
\end{figure}

Let $k$ be an non-negative integer. We use $F_k$ to denote the graph on $2k+1$ vertices consisting of $k$ triangles all sharing exactly one same vertex (an $F_k$ is also known as a {\it friendship graph}). Note that $F_0=K_1$.
Let $F_k^+$ be the graph on $4k+2$ vertices obtained by adding one pendant edge to each vertex in $F_k$. In other words, 
$$V(F_k^+)=\{a,a'\}\cup\bigcup_{i=1}^{k}\{b_i,c_i,b_i',c_i'\}\mbox{,}$$  
$$E(F_k^+)=\{aa'\}\cup\bigcup_{i=1}^{k}\{ab_i,ac_i,b_ic_i,b_ib_i',c_ic_i'\}\mbox{.}$$
For any positive integer $n$, there exists a unique pair of $r\in[0,3]$ and an non-negative integer $k$, such that $n=4k+r$. We define the graph $Sat_n$ on $n=4k+r$ vertices as follows. If $r=0$, $Sat_{n}=F_{k}^{+}-\{a', b_1'\}$; if $r=1$, $Sat_{n}=F_{k}^{+}-\{a'\}$; if $r=2$, $Sat_{n}=F_{k}^{+}$; if $r=3$, $Sat_{n}=F_{k+1}^{+}-\{a', b_1', c_1'\}$. Note that, $|E(Sat_n)|=\lceil\frac{5n}{4}-\frac{3}{2}\rceil$. See Figure~\ref{sat4}.

\begin{obs}\label{4cons}
Let $n$ be a positive integer and $I\subseteq[3,+\infty)$ be an integer set. If $\{3, 4\}\cap I=\{4\}$ and $\{5, 6, 7\}\cap I\neq\emptyset$, then $Sat_n$ is a $\mathcal{C}_{I}$-saturated graph on $n$ vertices.
In particular, $\sat(n,\mathcal{C}_{I})\le\lceil\frac{5n}{4}-\frac{3}{2}\rceil$. 
\end{obs}
\begin{proof}
Firstly, since any block of $Sat_n$ is either $K_2$ or $C_3$, there is no copy of cycle in $\mathcal{C}_{I}$ contained in $Sat_n$.
Also, one can check that, if two vertices $x,y\in V(Sat_n)$ are not adjacent, then either there is a path of length $3$ from $x$ to $y$, or we can suppose $\{x,y\}=\{b_i',b_j'\}$ with $i\neq j$ without loss of generality. In the first case, the addition of $xy$ would create a copy of $C_4\in\mathcal{C}_I$. In the second case, the addition of $b_i'b_j'$ would create the cycles $b_i'b_iab_jb_j'b_i'\cong C_5$, $b_i'b_ic_iab_jb_j'b_i'\cong C_6$ and $b_i'b_ic_iac_jb_jb_j'b_i'\cong C_7$. This completes the proof.
\end{proof}
Apparently, by Observation~\ref{4cons}, $\sat(n,\mathcal{C}_{\{4,5\}})\le\lceil\frac{5n}{4}-\frac{3}{2}\rceil$.

\subsection{Lower bound}
For any positive integer $i$ and any graph $H$, let $D_i(H)$ be the set of all the vertices of degree $i$ in $H$. We call any vertex in $D_1(H)$ a leaf in $H$.

In the following of this section, let $G$ be the smallest counterexample of Theorem~\ref{45}. 
To be more specific, suppose $G$ is a $\mathcal{C}_{\{4,5\}}$-saturated graph on $n\ge 1$ vertices with $|E(G)|<\lceil\frac{5n}{4}-\frac{3}{2}\rceil$ but $\sat(m,\mathcal{C}_{\{4,5\}})=\lceil\frac{5m}{4}-\frac{3}{2}\rceil$ for any $m\in[1,n-1]$. We suppose $n\ge 4$ since these small cases are easy to check.

\subsubsection{The vertices of degree $1$}
A {\it matching} is a graph whose components are all $K_2$'s or $K_1$'s. A {\it perfect matching} is a graph whose components are all $K_2$'s. 
Since $G$ does not contain a copy of $C_4$ as its subgraph, there is no copy of $P_3$ in the neighborhood of any vertex in $G$, so we have the following fact.
\begin{fact}\label{fact1}
for any vetex $v\in V(G)$, $G[N_{G}(v)]$ is a matching.
\end{fact}
The next three facts are also trivial by the saturation of $G$.
\begin{fact}\label{fact0}
$G$ is connected.
\end{fact}
\begin{fact}\label{fact2}
Suppose $u,v\in V(G)$ and $u\neq v$. If $uv\notin E(G)$, then there is a $P_4$ or a $P_5$ from $u$ to $v$. 
\end{fact}
\begin{fact}\label{fact3}
Suppose $u,v\in D_1(G)$ and $u\neq v$. Then $uv\notin E(G)$ and $N_{G}(u)\neq N_{G}(v)$. In other words, any vertex in $G$ has at most one neighbor of degree $1$.
\end{fact}
Now we give a lemma that characterizes the neighborhood of the only neighbor of any leaf in $G$.
\begin{lem}\label{1matching}
Let $u_0\in D_1(G)$ and $N_{G}(u_0)=\{v_0\}$. It holds that $G[N_{G}(v_0)\backslash\{u_0\}]$ is a perfect matching.
\end{lem}
\begin{proof}
Suppose $v$ is a neighbor of $v_0$ other than $u_0$. Apparently, $u_0v\notin E(G)$, so there is a path of length $3$ or $4$ from $u_0$ to $v$ by Fact~\ref{fact2}. Since the only neighbor of $u_0$ is $v_0$, there must be a path from $v_0$ to $v$ of length $2$ or $3$.
If there is a path $P$ of length $3$ from $v_0$ to $v$, then $P\cup\{v_0v\}$ is a copy of $C_4$ in $G$, a contradiction by the saturation of $G$. Hence, there must be a path of length $2$ from $v_0$ to $v$, so $v_0$ and $v$ have a common neighbor, say $w$. Now $vw\in E(G[N_{G}(v_0)\backslash\{u_0\}])$, which implies that any vertex $v\in N_{G}(v_0)\backslash\{u_0\}$ is not contained in a component isomorphic to $K_1$ in $G[N_{G}(v_0)\backslash\{u_0\}]$. By Fact~\ref{fact1}, we are done.  
\end{proof}

For a connected graph $H$ and two vertices $x,y\in V(H)$, the {\it distance} of $x$ and $y$, denoted by $d_{H}(x,y)$, is the smallest possible length of a path in $H$ from $x$ to $y$. 
For any vertex $u\in V(H)$ and non-negative integer $i$, we use $L_{i}^H(u)$ to denote the set consisting of all the vertices having distance $i$ with $u$. Note that $L_{0}^{H}(u)=\{u\}$ and $L_{1}^{H}(u)=N_H(u)$.
For any vertex $u\in V(H)$, a {\it BFS-tree} (or a breadth-first search tree) of $H$ rooted at $u$ is a spanning tree $T$ of $H$ such that $d_T(u,v)=d_H(u,v)$ for any $v\in V(H)$. It is easy to see that $T[L_i^H(u)]$ is an empty graph for any non-negative integer $i$ if $T$ is a BFS-tree of $H$ rooted at $u$.

The following lemma gives a bound of the degree of the only neighbor of a leaf in $G$.
\begin{lem}\label{1degree}
Let $u_0\in D_1(G)$ and $N_{G}(u_0)=\{v_0\}$. It holds $d_{G}(v_0)\ge 5$.
\end{lem}
\begin{proof}
Suppose this does not hold. Then by Fact~\ref{fact3} and Lemma~\ref{1matching}, it must hold that $d_{G}(v_0)=3$ and $G[N_{G}(v_0)\backslash\{u_0\}]\cong K_2$.
Suppose $N_{G}(v_0)=\{u_0,x,y\}$ and $xy\in E(G)$. Note that $L_0^G(u_0)=\{u_0\}$, $L_1^G(u_0)=\{v_0\}$ and $L_2^G(u_0)=\{x,y\}$. Hence, every vertex in $L_3^G(u_0)$ is a neighbor of $x$ or $y$ in $G$. Let $X=N_{G}(x)\cap L_3^G(u_0)$ and $Y=N_{G}(y)\cap L_3^G(u_0)$.

Let $A_i=D_i(G[L_3^G(u_0)])$ for $i=0,1$. Let $B_1=D_1(G)\cap L_4^G(u_0)$ and $B_{\ge 2}=L_4^G(u_0)\backslash B_1$.
Now we give some claims.
\begin{claim}\label{c1}
for any $i\ge 5$, $L_i^{G}(u_0)=\emptyset$.
\end{claim}
This claim can be directly deduced by Fact~\ref{fact2}.
\begin{claim}\label{c2}
$X\cap Y=\emptyset$ and $G[X,Y]$ is an empty graph. Moreover, $G[L_3^G(u_0)]$ is a matching, and so $L_3^G(u_0)=A_0\cup A_1$.
\end{claim}
It is easy to see that $X\cap Y=\emptyset$.
Otherwise, suppose $z\in X\cap Y$, then $v_0xzyv_0\cong C_4$ is contained in $G$, a contradiction.
Also, $G[X,Y]$ is an empty graph. Otherwise, if there is an edge $x'y'\in E(G)$ with $x'\in X$ and $y'\in Y$, then $v_0xx'y'yv_0\cong C_5$ is contained in $G$, a contradiction, too. 
Thus, by Fact~\ref{fact1}, $G[L_3^G(u_0)]=G[X]\cup G[Y]$ is a matching. The claim is done.
\begin{claim}\label{c3}
The only neighbor of any vertex $v\in B_1$ must be in $A_1$.
\end{claim}
Suppose $N_{G}(v)=\{w\}$. By defintion, every vertex in $L_4^G(u_0)$ has at least one neighbor in $L_3^G(u_0)=A_0\cup A_1$, so $w\in A_0\cup A_1$. 
Since $X\cap Y=\emptyset$ by Claim~\ref{c2}, suppose $w\in X$ and $w\notin Y$ without loss of generality. Thus, $x\in N_{G}(w)\backslash\{v\}$. 
If $w\in A_0$,
then $w$ has no neighbor in $L_3^G(u_0)\cup\{y,v_0\}$, which implies that $w$ and $x$ has no common neighbor. Hence, $x$ is an isolated vertex in $G[N_{G}(w)\backslash\{v\}]$. By Lemma~\ref{1matching}, this is a contradiction. So $w\in A_1$ and the claim is done.
\begin{claim}\label{c4}
If $a, b\in L_3^G(u_0)$ and $a\neq b$, then $a$ and $b$ do not have common neighbors in $L_4^G(u_0)$.
\end{claim}
Suppose to the contrary that $c\in N_{G}(a)\cap N_{G}(b)\cap L_4^G(u_0)$. Without loss of generality, suppose $a\in X$. If $b\in X$, then $xacbx\cong C_4$ is contained in $G$, a contradiction. If $b\in Y$, then $xacbyx\cong C_5$ is contained in $G$, a contradiction, too. The claim is done.
\begin{claim}\label{c5}
$|A_1|+|B_{\ge 2}|\ge\frac{n}{2}-3$.
\end{claim}
For any vertex $v\in A_0$, if $N_{G}(v)\cap L_4^G(u_0)=\emptyset$, then $N_{G}(v)\subset L_2^G(u_0)=\{x,y\}$. By Claim~\ref{c2}, this means $d_G(v)=1$.
Since $A_0\subseteq L_3^G(u_0)\subset N_{G}(x)\cup N_{G}(y)$, by Fact~\ref{fact3}, there are at most two such vertices $v\in A_0$ with $N_{G}(v)\cap L_4^G(u_0)=\emptyset$.
Hence, by Claim~\ref{c3} and Claim~\ref{c4}, it holds $|B_{\ge 2}|\ge|N_{G}(A_0)\cap L_4^G(u_0)|\ge |A_0|-2$.
Also, by Fact~\ref{fact3} and Claim~\ref{c3}, it holds $|A_1|\ge |B_1|$. Put them together, we have
$$|A_1|+|B_{\ge 2}|\ge|A_0|+|B_1|-2\mbox{.}$$
On the other hand, 
$$(|A_1|+|B_{\ge 2}|)+(|A_0|+|B_1|-2)=|V(G)\backslash\{u_0,v_0,x,y\}|-2=n-6\mbox{.}$$
Therefore,
$$|A_1|+|B_{\ge 2}|\ge\frac{1}{2}((|A_1|+|B_{\ge 2}|)+(|A_0|+|B_1|-2))=\frac{n}{2}-3\mbox{.}$$
The claim is done.

Now we consider a BFS-tree $T$ rooted at $u_0$ in $G$. Since $T[L_i^G(u_0)]$ is an empty graph for $i\in[2,3]$ and the vertices in $L_4^G(u_0)$ are all leaves in $T$ by Claim~\ref{c1}, we have
\begin{equation*}
\begin{split}
|E(G)| &\ge |E(T)|+|E(G[L_2^G(u_0)])|+|E(G[L_3^G(u_0)])|+\frac{1}{2}\sum_{v\in L_4^G(u_0)}(d_{G}(v)-1)\\
       &\ge (n-1)+1+\frac{1}{2}|A_1|+\frac{1}{2}|B_{\ge 2}|\\
       &\ge \frac{5n}{4}-\frac{3}{2}\mbox{,}
\end{split}
\end{equation*}
where the last inequality holds by Claim~\ref{c5}. By the assumption that $|E(G)|<\lceil\frac{5n}{4}-\frac{3}{2}\rceil$ and $|E(G)|$ being an integer, this is a contradiction.
\end{proof}

\subsubsection{The vertices of degree $2$}
We call a path $P$ in $G$ a {\it degenerated path} if $V(P)\subseteq D_2(G)$. 
For any degenerated path $P=x_1x_2\dots x_r$, since $d_{G}(x_i)=2$ for every $i\in[1,r]$, it holds $|N_G(x_j)\backslash V(P)|=1$ for $j=1,r$. Let $N_G(x_1)\backslash V(P)=\{x_0\}$ and $N_G(x_r)\backslash V(P)=\{x_{r+1}\}$. We use $A(P)$ to denote the path (if $x_0\neq x_{r+1}$) or the cycle (if $x_0= x_{r+1}$) $x_0x_1\dots x_rx_{r+1}$. The following fact is trivial.

\begin{fact}\label{fact5}
Let $P$ be a degenerated path in $G$.\\
(i) For any cycle $C$ in $G$, if $V(C)\cap V(P)\neq\emptyset$, then $A(P)\subseteq C$.\\
(ii) Let $ab$ be a non-edge in $G$ with $\{a,b\}\cap V(P)=\emptyset$. For any cycle $C$ in $G+ab$, if $V(C)\cap V(P)\neq\emptyset$, then $A(P)\subseteq C$. In other words, For any path $P'$ from $a$ to $b$ in $G$, if $V(P')\cap V(P)\neq\emptyset$, then $A(P)\subseteq P'$.
\end{fact}
Here are some lemmas for degenerated paths in $G$.

\begin{lem}\label{p<4}
Every degenerated path in $G$ has length at most $2$.
\end{lem}
\begin{proof}
Suppose to the contrary that there is a degenerated path $P=x_1x_2x_3x_4$ of length $3$. Let $A(P)=x_0x_1x_2x_3x_4x_5$.
Consider $G'=G-V(P)$. By defintion, the addition of any non-edge $e$ in $G'$ would create a cycle $C$ of length $4$ or $5$ in $G+e$. If $V(C)\cap V(P)\neq \emptyset$, by Fact~\ref{fact5}, $A(P)\subseteq V(C)$. If $x_0=x_5$, then $C_5\cong A(P)\subset G$, a contradiction by the saturation of $G$. So $x_0\neq x_5$. However, this means $|V(C)|\ge|V(A(P))|=6$, a contradiction by the length of $C$ being $4$ or $5$. 
Therefore, $V(C)\cap V(P)=\emptyset$. This implies that $G'$ is also a $\mathcal{C}_{\{4,5\}}$-saturated graph. However,
$$|E(G')|=|E(G)|-5<\lceil\frac{5}{4}n-\frac{3}{2}\rceil-5= \lceil\frac{5}{4}(n-4)-\frac{3}{2}\rceil=\lceil\frac{5}{4}|V(G')|-\frac{3}{2}\rceil\mbox{,}$$
which is a contradiction by the minimality of $G$. 
\end{proof}

\begin{lem}\label{c3block}
Every degenerated path $P$ of length $1$ is not contained in a copy of $C_3$ in $G$.
\end{lem}
\begin{proof}
Let $P=x_1x_2$. Note that $A(P)$ is a triangle, so we can suppose $A(P)=x_0x_1x_2x_0$. Hence, $x_0$ is a cut vertex of $G$ and the triangle $x_0x_1x_2x_0$ is a block. In this case, let $G'=G-\{x_1,x_2\}$. By definiton, the addition of any non-edge $e$ in $G'$ would create a cycle $C$ of length $4$ or $5$ in $G+e$ and apparently $x_1,x_2\notin V(C)$ (or $C$ would contain a triangle as its subgraph, which is impossible). This implies that $G'$ is also $\mathcal{C}_{\{4,5\}}$-saturated.
However, 
$$|E(G')|=|E(G)|-3<\lceil\frac{5}{4}n-\frac{3}{2}\rceil-3\le \lceil\frac{5}{4}(n-2)-\frac{3}{2}\rceil=\lceil\frac{5}{4}|V(G')|-\frac{3}{2}\rceil\mbox{,}$$
which is a contradiction by the minimality of $G$.
\end{proof}

\begin{lem}\label{pc6}
Every degenerated path $P$ of length at least $1$ must be contained in an induced copy $C$ of $C_6$ in $G$. In particular, $A(P)\subset C$ is a path.
\end{lem}
\begin{proof}
Let $P=x_1x_2\dots x_r$ be a degenerated path of length at least $1$ in $G$. Let $A(P)=x_0x_1\dots x_rx_{r+1}$.
Since $r\ge 2$, it holds $d_G(x_1)=d_{G}(x_2)=2$. Consider the pair $\{x_0,x_2\}$.
If $x_0x_2\in E(G)$, then $r=2$ and $x_3=x_0$. Hence, the degenerated path $x_1x_2$ is contained in a triangle $x_0x_1x_2x_0$, a contradiction by Lemma~\ref{c3block}. \\
Hence, $x_0x_2\notin E(G)$. By Fact~\ref{fact2}, there is a path $P'$ of length $3$ or $4$ from $x_0$ to $x_2$. Since $x_1\notin V(P')$ (or the length of $P'$ would be $2$), $x_0x_1x_2\cup P'$ is a cycle of length $2+3=5$ or $2+4=6$. Since $G$ cannot contain a copy of $C_5$ as its subgraph, $x_0x_1x_2\cup P'\cong C_6$. Thus, $x_0x_1x_2$ is in a copy of $C_6$ in $G$, which then implies that $P$ must also be in a copy of $C_6$ in $G$ by Fact~\ref{fact5}. Apparently, this copy of $C_6$ has no chord (or we can find copies of $C_4$ or $C_5$ within its vertices), so it must be induced.
\end{proof}
By Lemma~\ref{pc6}, for any degenerated path $P$ of length at least $1$, there are some induced cycles of length $6$ containing $A(P)$. We use $C_6(P)$ to denote one of them. 

\begin{figure}[h]
\centering
\subfigure[if $P$ has length $2$]{\includegraphics[width=1.96in]{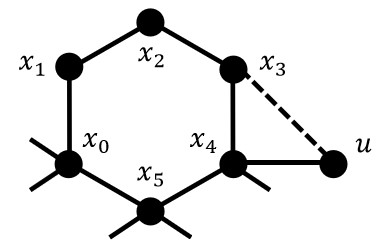}}
\hspace{0.8in}
\subfigure[if $P$ has length $1$]{\includegraphics[width=1.88in]{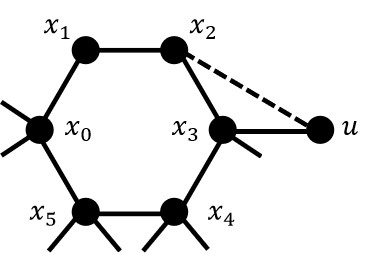}}
\caption{The proof of Lemma~\ref{apc3}}\label{proof2}
\end{figure}

\begin{lem}\label{apc3}
Let $P$ be a degenerated path of length $1$ or $2$ in $G$ and pick $v\in V(A(P))\backslash V(P)$. It holds that either $G[N_{G}(v)\backslash V(C_6(P))]$ or $G[N_{G}(v)\backslash V(P)]$ is a perfect matching.
\end{lem}
\begin{proof}
By Lemma~\ref{pc6}, $N_{G}(v)\backslash V(C_6(P))$ can be obtained by deleting one vertex from $N_{G}(v)\backslash V(P)$.
By Fact~\ref{fact1}, for any $v\in V(A(P))\backslash V(P)$, to prove either $G[N_{G}(v)\backslash V(C_6(P))]$ or $G[N_{G}(v)\backslash V(P)]$ is a perfect matching, we only need to prove that every vertex in $N_{G}(v) \backslash V(C_6(P))$ has a neighbor in $N_{G}(v)\backslash V(P)$.\\
We consider the case that $P$ has length $2$ at first. Let $P=x_1x_2x_3$, $A(P)=x_0x_1x_2x_3x_4$ and $C_6(P)=x_5x_0x_1x_2x_3x_4x_5$. We pick $x_4\in V(A(P))\backslash V(P)$ without loss of generality. Let $u\in N_{G}(x_4)\backslash V(C_6(P))$.
Apparently, $ux_3$ is a non-edge in $G$. By Fact~\ref{fact2}, there is a path $Q$ of length $3$ or $4$ from $u$ to $x_3$ in $G$. Note that $x_1x_2$ is a degenerated path of length $1$. By Fact~\ref{fact5}, if $V(Q)\cap V(P)\neq \emptyset$, then $\{u\}\cup V(A(x_1x_2))\subseteq V(Q)$.
Since $|V(Q)|\le 5=|\{u\}\cup V(A(x_1x_2))|$, we must have $V(Q)=\{u\}\cup V(A(x_1x_2))$, which implies $Q=ux_0x_1x_2x_3$. Thus $u\in N_G(x_0)$, which
gives a cycle $ux_0x_5x_4u$ of length $4$ in $G$, a contradiction. So $V(Q)\cap V(P)=\emptyset$. In particular, $x_2\notin V(Q)$. Since $x_3$ is an endpoint of $Q$ and $x_4$ is the only neighbor of $x_3$ other than $x_2$, we have $x_3x_4\in E(Q)$. Thus, there is a path $Q'$ of length $2$ or $3$ from $u$ to $x_4$, where $Q'=Q-\{x_3\}$.
If the length of $Q'$ is $3$, then $ux_4\cup Q'$ is a copy of $C_4$ in $G$, a contradiction. So $Q'$ has length $2$, which means $u$ has a neighbor in $N_{G}(x_4)\backslash V(P)$.\\
Now consider the case that $P$ has length $1$.  Let $P=x_1x_2$, $A(P)=x_0x_1x_2x_3$ and $C_6(P)=x_4x_5x_0x_1x_2x_3x_4$. We pick $x_3\in V(A(P))\backslash V(P)$ without loss of generality. Let $u\in N_{G}(x_3)\backslash V(C_6(P))$. Apparently, $ux_2$ is a non-edge in $G$. By Fact~\ref{fact2}, there is a path $Q$ of length $3$ or $4$ from $u$ to $x_2$ in $G$. \\
Suppose $V(Q)\cap V(P)\neq \emptyset$.
Since $x_1$ is a degenerated path of length $0$, by Fact~\ref{fact5}, this means $\{u\}\cup V(A(P))\subseteq V(Q)$. Note that $|V(Q)|-1\le|\{u\}\cup V(A(P))|=4\le|V(Q)|$. If $|V(Q)|=|\{u\}\cup V(A(P))|=4$, then $Q=ux_0x_1x_2$. Thus $u\in N_G(x_0)$, which gives a cycle $ux_0x_5x_4x_3u$ of length $5$ in $G$, a contradiction. So $|V(Q)|=5$ and we can suppose $Q=uwx_0x_1x_2$, where $w\notin V(C_6(P))\cup\{u\}$. In this case, consider $G'=G-V(P)=G-\{x_1,x_2\}$. By definiton, the addition of any non-edge $ab$ in $G'$ would create a cycle $C$ of length $4$ or $5$ in $G+ab$.
Claim that: we can always pick a cycle $C'$ in $G'+ab$ of length $4$ or $5$.
If $V(P)\cap V(C)=\emptyset$, pick $C'=C$ and the claim is done.
If $V(P)\cap V(C)\neq\emptyset$, then by Fact~\ref{fact5}, $\{a,b\}\cup V(A(P))\subseteq V(C)$. So $|\{a,b\}\cup V(A(P))|\le |V(C)|\le 5=|\{a,b\}| +  |V(A(P))|-1$, which implies $\{a,b\}\cap V(A(P))\neq\emptyset$. Note that $a,b\in V(G)\backslash\{x_1,x_2\}$. 
If $|\{a,b\}\cap V(A(P))|=2$, then $\{a,b\}=\{x_0, x_3\}$. In this case, we pick $C'=x_0x_3x_4x_5x_0\cong C_4$.
If $|\{a,b\}\cap V(A(P))|=1$, we can suppose $a=x_0$ without loss of generality (the other cases are similar). 
There is at least one of the two sets $\{x_4, x_5\}$ and $\{u, w\}$ that does not contain $b$, suppose it is $\{x_4, x_5\}$ without loss of generality. Then pick $C'=(C\backslash x_0x_1x_2x_3 )\cup x_0x_5x_4x_3\cong C$. Hence, the claim is done. This then implies that $G'$ is also a $\mathcal{C}_{\{4,5\}}$-saturated graph. However,
$$|E(G')|=|E(G)|-3<\lceil\frac{5}{4}n-\frac{3}{2}\rceil-3\le \lceil\frac{5}{4}(n-2)-\frac{3}{2}\rceil=\lceil\frac{5}{4}|V(G')|-\frac{3}{2}\rceil\mbox{,}$$
which is a contradiction by the minimality of $G$.\\
Therefore, $V(Q)\cap V(P)=\emptyset$. Similarly as the previous proof, there is a path $Q'$ of length $2$ or $3$ from $u$ to $x_3$, where $Q'=Q-\{x_2\}$. If the length of $Q'$ is $3$, then $ux_3\cup Q'$ is a copy of $C_4$ in $G$, a contradiction. So $Q'$ has length $2$, which means $u$ has a neighbor in $N_{G}(x_3)\backslash V(P)$. See Figure~\ref{proof2}.
\end{proof}

For $i\in[0,2]$, let $D_2^i(G)=\{v\in D_2(G): |N_G(v)\cap D_2(G)|=i\}$. 
Let $D_2^{1+}(G)=\{v\in D_2^1(G): N_{G}(v)\cap D_2^2(G)\neq\emptyset\}$ and $D_2^{1-}(G)=D_2^1(G)\backslash D_2^{1+}(G)$.
By Lemma~\ref{p<4}, the only neighbor of degree $2$ of every vertex in $D_2^{1+}(G)$ comes from $D_2^2(G)$, the only neighbor of degree $2$ of every vertex in $D_2^{1-}(G)$ comes from $D_2^{1-}(G)$, and the two neighbors of every vertex in $D_2^2(G)$ is in $D_2^{1+}(G)$. 
The next lemma is about a special case for vertices in $D_2^0(G)$.

\begin{figure}[h]
\centering
\includegraphics[width=1.7in]{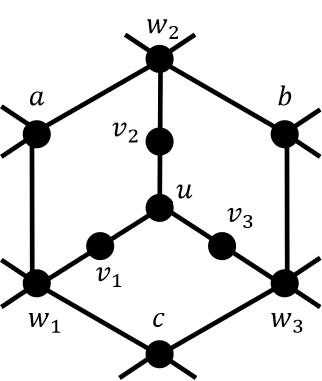}
\caption{The proof of Lemma~\ref{3222}}\label{proof23}
\end{figure}

\begin{lem}\label{3222}
Let $u$ be a vertex of degree $3$ in $G$ and $N_{G}(u)=\{v_1,v_2,v_3\}\subseteq D_2^0(G)$. Let $N_{G}(v_i)=\{u,w_i\}$ for $i\in[1,3]$. Then $w_1, w_2, w_3$ are pairwise distinct, pairwise non-adjacent and $N_{G}(w_1)\cap N_{G}(w_2)\cap N_{G}(w_3)\neq\emptyset$. 
\end{lem}
\begin{proof}
If $w_1, w_2, w_3$ are not pairwise distinct, we can suppose $w_1=w_2$ without loss of generality. Then $w_1v_1uv_2w_1$ is a cycle of length $4$ in $G$, a contradicition. 
Similarly, if $w_1, w_2, w_3$ are not pairwise non-adjacent, we can suppose $w_1w_2\in E(G)$ without loss of generality. Then $w_1v_1uv_2w_2w_1$ is a cycle of length $5$ in $G$, a contradicition.
So it remains to prove $N_{G}(w_1)\cap N_{G}(w_2)\cap N_{G}(w_3)\neq\emptyset$.\\
Consider the non-edge $v_1v_2$. By Fact~\ref{fact2}, there is a path $P$ of length $3$ or $4$ from $v_1$ to $v_2$ in $G$. The vertex adjacent to $v_1$ within $P$ must be in $N_{G}(v_1)=\{u,w_1\}$. In other words, either $v_1u\in E(P)$ or $v_1w_1\in E(P)$.
If $v_1u\in E(P)$, since $v_3$ is the only neighbor of $u$ other than $v_1$ and $v_2$, it holds $uv_3\in E(P)$. Then in a similar way, $v_3w_3\in E(P)$. Since $|E(P)|\le 4$ and $v_2\in V(P)$, we then must have $w_3v_2\in E(P)$. However, this means $w_2=w_3$, a contradiction by the previous proof.
Hence, $v_1w_1\in E(P)$. Similarly, $w_2v_2\in E(P)$. If $w_1w_2\in E(G)$, then $uv_1w_1w_2v_2u$ is a cycle of length $5$ in $G$, a contradiction. So $d_{G}(w_1,w_2)\ge 2$. Since $|E(P)|\le 4$, we must have $d_{G}(w_1,w_2)= 2$ and $P=v_1w_1aw_2v_2$ where $a\in N_{G}(w_1)\cap N_{G}(w_2)$. Similar, we can pick $b\in N_{G}(w_2)\cap N_{G}(w_3)$ and $c\in N_{G}(w_3)\cap N_{G}(w_1)$. \\
If $a,b,c$ are not pairwise distinct, suppose $a=b$ without loss of generality. We then have $a=b\in N_{G}(w_1)\cap N_{G}(w_2)\cap N_{G}(w_3)$ and we are done. So $a, b, c$ are pairwise distinct and $w_1aw_2bw_3cw_1$ is a cycle of length $6$ in $G$.
Let $G'=G-\{u, v_1, v_2, v_3\}$. By definiton, the addition of any non-edge $e$ in $G'$ would create a cycle $C$ of length $4$ or $5$ in $G+e$. 
If $V(C)\cap \{u, v_1, v_2, v_3\}\neq\emptyset$, it is easy to check that $V(C)$ contains at least two members of $\{v_1,v_2,v_3\}$. Without loss of generality, suppose $v_1,v_2\in V(C)$. By Fact~\ref{fact5}, $\{w_1,v_1,u,v_2,w_2\}=V(A(v_1)\cup A(v_2))\subseteq V(C)$. Since $|V(C)|\le 5$, we have $V(C)=\{w_1,v_1,u,v_2,w_2\}$ and $C=w_1v_1uv_2w_2w_1$, which leads to $e=w_1w_2$. In this case, we can find $w_1w_2bw_3cw_1\cong C_5$ in $G'+e$.
This implies that $G'$ is also $\mathcal{C}_{\{4,5\}}$-saturated.
However, 
$$|E(G')|=|E(G)|-6<\lceil\frac{5}{4}n-\frac{3}{2}\rceil-6< \lceil\frac{5}{4}(n-4)-\frac{3}{2}\rceil=\lceil\frac{5}{4}|V(G')|-\frac{3}{2}\rceil\mbox{,}$$
which is a contradiction by the minimality of $G$. See Figure~\ref{proof23}.
\end{proof}

\subsubsection{The discharging method}
\begin{figure}[h]
\centering
\includegraphics[width=5.5in]{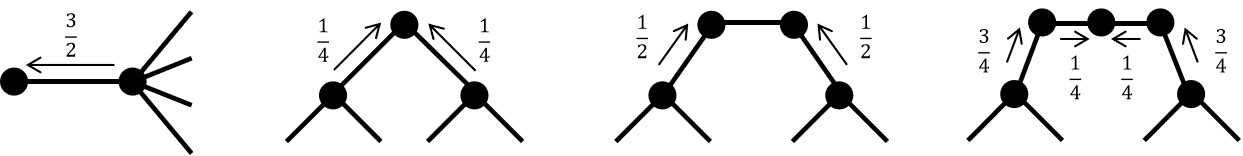}
\caption{The rules}\label{charge}
\end{figure}

For each vertex $v\in V(G)$, define its initial charge as $ch(v)=d_{G}(v)-\frac{5}{2}$. Then
$$\sum_{v\in V(G)}ch(v)=\sum_{v\in V(G)}d_{G}(v)-\frac{5}{2}|V(G)|=2|E(G)|-\frac{5n}{2}\mbox{.}$$
To get a contradiction by $|E(G)|<\lceil\frac{5}{4}n-\frac{3}{2}\rceil$, we will prove $\sum_{v\in V(G)}ch(v)\ge -\frac{1}{4}$. Now we redistribute the charges according to the following rules. See Figure~\ref{charge}.

(R1) Every vertex $v\in D_{2}^0(G)\cup D_{2}^2(G)$ gets $\frac{1}{4}$ from each of its two neighbors.

(R2) Every vertex $v\in D_{2}^{1+}(G)$ gets $\frac{3}{4}$ from its only neighbor of degree more than $2$.

(R3) Every vertex $v\in D_{2}^{1-}(G)$ gets $\frac{1}{2}$ from its only neighbor of degree more than $2$.

(R4) Every vertex $v\in D_1(G)$ gets $\frac{3}{2}$ from its only neighbor.\\
Let the charge of every vertex $v\in V(G)$ end up with $ch'(v)$. We have the following observation and lemmas.

\begin{obs}\label{ch12}
For vertex $v\in D_1(G)\cup D_2(G)$, $ch'(v)=0$.
\end{obs}
\begin{proof}
By the rules, if $v\in D_{2}^0(G)\cup D_{2}^2(G)$, then $ch'(v)=ch(v)+2\times\frac{1}{4}=0$; if $v\in D_{2}^{1+}(G)$, then $ch'(v)=ch(v)+\frac{3}{4}-\frac{1}{4}=0$; if $v\in D_{2}^{1-}(G)$, then $ch'(v)=ch(v)+\frac{1}{2}=0$; if $v\in D_{1}(G)$, then $ch'(v)=ch(v)+\frac{3}{2}=0$.
\end{proof}

\begin{lem}\label{ch1n}
For any vertex $v\in V(G)\backslash (D_1(G)\cup D_2(G))$ with $N_{G}(v)\cap D_1(G)\neq \emptyset$, $ch'(v)\ge\frac{1}{2}$.
\end{lem}
\begin{proof}
By Fact~\ref{fact3}, there is exact one vertex of degree $1$ in $N_{G}(v)$, so we can suppose $N_{G}(v)\cap D_1(G)=\{u\}$.
By Lemma~\ref{1matching}, $G[N_{G}(v)\backslash\{u\}]$ is a perfect matching.
Thus, for any $x\in D_2(G)\cap(N_{G}(v)\backslash\{u\})$, $N_{G}(x)=\{v,y\}$, where $y\in N_{G}(v)\backslash\{u\}$. So $xyvx$ is a copy of $C_3$.
By Lemma~\ref{c3block} and Fact~\ref{fact5}, $y\notin D_2(G)$ and $x\in D_2^0(G)$.
So $v$ should only give $\frac{3}{2}$ to $u$ and give $\frac{1}{4}$ to at most half of the vertices in $N_{G}(v)\backslash\{u\}$.
By Lemma~\ref{1degree}, 
$$ch'(v)\ge (d_{G}(v)-\frac{5}{2})-\frac{3}{2}-\frac{d_{G}(v)-1}{2}\times\frac{1}{4}=\frac{7d_{G}(v)-31}{8}\ge \frac{1}{2}\mbox{,}$$
which completes the proof.
\end{proof}
For any vertex $v\in V(G)\backslash (D_1(G)\cup D_2(G))$ with $N_{G}(v)\cap D_1(G)=\emptyset$. Let $r_G(v)=|N_{G}(v)\cap  D_2^0(G)|$, $s^+_G(v)=|N_{G}(v)\cap  D_2^{1+}(G)|$, $s^-_G(v)=|N_{G}(v)\cap  D_2^{1-}(G)|$ and $t_G(v)=|N_{G}(v)\backslash D_2(G)|$.

\begin{lem}\label{chother}
Let $v\in V(G)\backslash (D_1(G)\cup D_2(G))$ with $N_{G}(v)\cap D_1(G)=\emptyset$. Let $r=r_G(v)$, $s^+=s^+_G(v)$, $s^-=s^-_G(v)$ and $t=t_G(v)$, then:\\
(i) $d_G(v)=r+s^++s^-+t$ and $2s^++s^-\le 2$;\\
(ii) $ch'(v)=t+\frac{3}{4}r+\frac{1}{2}s^-+\frac{1}{4}s^+-\frac{5}{2}$.
\end{lem}
\begin{proof}
it is trivial to check that $d_G(v)=r+s^++s^-+t$ and $ch'(v)=t+\frac{3}{4}r+\frac{1}{2}s^-+\frac{1}{4}s^+-\frac{5}{2}$. 
It remains to prove $2s^++s^-\le 2$.\\
If $2s^++s^->0$, then $v$ is adjacent to some vertices in the degenerated path $P$ of length at least $1$.
By Lemma~\ref{apc3}, all the vertices in $N_{G}(v)\backslash V(C_6(P))$ is in a copy of $C_3$, which implies that $(N_{G}(v)\backslash V(C_6(P)))\cap D_2^1(G)=\emptyset$ by Lemma~\ref{c3block} and Fact~\ref{fact5}. 
Since $C_6(P)$ is an induced cycle by Lemma~\ref{pc6}, $N_{G}(v)\cap D_2^1(G)\subseteq N_{C_6(P)}(v)$.
Let $N_{C_6(P)}(v)=\{u_1, u_5\}$ and $C_6(P)=vu_1u_2u_3u_4u_5v$.
If $\{u_1, u_5\}\not\subseteq D_2^1(G)$, then $s^++s^-\le 1$ and the proof is done.
If $\{u_1, u_5\}\subseteq D_2^1(G)$, then $u_2, u_4\in D_2^{1-}(G)\cup D_2^2(G)$.
If one of $u_2$ and $u_4$ is in $D_2^2(G)$, then $u_3\in D_2(G)$, which means $u_1u_2u_3u_4u_5$ is a degenerated path of length $4$, a contradiction by Lemma~\ref{p<4}. Hence, $u_2, u_4\in D_2^{1-}(G)$, so $s^+=0$ and $s^-=2$. This completes the proof.
\end{proof}

\begin{figure}[h]
\centering
\subfigure[if $s^+\ge 1$]{\includegraphics[width=1.68in]{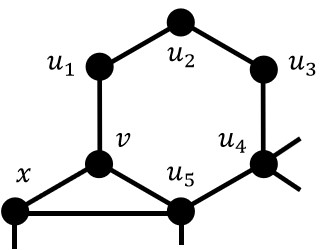}}
\hspace{0.8in}
\subfigure[if $s^+=0$ and $s^-\ge 1$]{\includegraphics[width=1.5in]{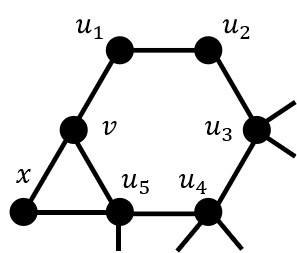}}
\caption{The proof of Lemma~\ref{ch4}}\label{proof31}
\end{figure}

\begin{lem}\label{ch4}
If $v\in V(G)$ and $ch'(v)<0$, then $v\in D_3(G)$. Moreover, $N_{G}(v)\subseteq D_2^0(G)$.
\end{lem}
\begin{proof}
We firstly prove that $v\in D_3(G)$. Suppose to the contrary that $v\notin D_3(G)$, by Observation~\ref{ch12}, $d_{G}(v)\ge 4$. 
Note that $N_{G}(v)\cap D_1(G)=\emptyset$ by Lemma~\ref{ch1n}. Let $r=r_G(v)$, $s^+=s^+_G(v)$, $s^-=s^-_G(v)$ and $t=t_G(v)$.
Then by Lemma~\ref{chother} (i) and (ii), 
\begin{equation*}
\begin{split}
ch'(v) &= t+\frac{3}{4}r+\frac{1}{2}s^-+\frac{1}{4}s^+-\frac{5}{2}\\
       &= \frac{1}{4}t+\frac{3}{4}(r+s^++s^-+t)-\frac{1}{4}(2s^++s^-)-\frac{5}{2}\\
       &\ge 0+\frac{3}{4}d_{G}(v)-\frac{1}{2}-\frac{5}{2}\\
       &= \frac{3}{4}d_{G}(v)-3\\
       &\ge 0\mbox{,}
\end{split}
\end{equation*}
a contradiction.\\
For the ``moreover" part,  suppose to the contrary that $N_{G}(v)\not\subseteq D_2^0(G)$. 
Since $v\in D_3(G)$, it holds $r\le 2$, $t+s^++s^-\ge 1$.
If $s^+=s^-=0$, 
by Lemma~\ref{chother} (ii), $ch'(v)=t+\frac{3}{4}r-\frac{5}{2}\ge\frac{3}{4}d_{G}(v)-\frac{5}{2}>0$, a contradiction.\\
If $s^+\ge 1$, then $v$ is adjacent to some vertices in the degenerated path $P$ of length $2$.
Let $P=u_1u_2u_3$ and $C_6(P)=vu_1u_2u_3u_4u_5v$. Since $v\in D_3(G)$, $N_{G}(v)=\{u_1, u_5, x\}$, where $x\notin V(C_6(P))$.
By Lemma~\ref{apc3}, we must have $u_5x\in E(G)$.
Now consider $G'=G-\{u_1, u_2, u_3\}$. By definition, the addition of any non-edge $ab$ in $G'$ would create a cycle $C$ of length $4$ or $5$ in $G+ab$.
Claim that: we can always pick a cycle $C'$ in $G'+ab$ of length $4$ or $5$.
If $V(P)\cap V(C)=\emptyset$, pick $C'=C$ and the claim is done. If $V(P)\cap V(C)\neq\emptyset$, by Fact~\ref{fact5}, $\{a,b\}\cup V(A(P))\subseteq V(C)$. Since $|V(A(P))|=5\ge |V(C)|$, we have $a,b\in V(A(P))$. Thus, $\{a,b\}=\{v, u_4\}$. In this case, we pick $C'=vu_4u_5xv$, which is a cycle of length $4$.
Hence, the claim is done.
This implies that $G'$ is also a $\mathcal{C}_{\{4,5\}}$-saturated graph. However,
$$|E(G')|=|E(G)|-4<\lceil\frac{5}{4}n-\frac{3}{2}\rceil-4\le \lceil\frac{5}{4}(n-3)-\frac{3}{2}\rceil=\lceil\frac{5}{4}|V(G')|-\frac{3}{2}\rceil\mbox{,}$$
which is a contradiction by the minimality of $G$.\\
Therefore, $s^+=0$ and $s^-\ge 1$, which means $v$ is adjacent to some vertices in the degenerated path $P$ of length $1$.
Let $P=u_1u_2$ and $C_6(P)=vu_1u_2u_3u_4u_5v$. Since $v\in D_3(G)$, $N_{G}(v)=\{u_1, u_5, x\}$, where $x\notin V(C_6(P))$.
By Lemma~\ref{apc3}, we must have $u_5x\in E(G)$.
Note that $d_{G}(x)\ge |\{u_5,v\}|=2$ and $d_{G}(u_5)\ge|\{u_4,v,x\}|=3$. If $d_{G}(x)\ge 3$, then $t=2$, $s^-=1$ and $r=s^+=0$, which implies $ch'(v)=0$ by Lemma~\ref{chother} (ii). This is a contradiction. So $d_{G}(x)=2$ and $N_{G}(x)=\{u_5,v\}$.\\
Now consider $G_1=G-\{x,v,u_1,u_2\}$ and $G_2=G_1+u_3u_5$. If $G_2$ does not contain cycles of length $4$ or $5$ put $G'=G_2$. Otherwise, let $G'=G_1$.
By definition, the addition of any non-edge $ab$ in $G'$ would create a cycle $C$ of length $4$ or $5$ in $G+ab$.
Note that $\{u_3,u_5\}$ is a cut set between $\{x,v,u_1,u_2\}$ and $V(G')\backslash \{u_3, u_5\}$. Since $\{a,b\}\cap\{x,v,u_1,u_2\}=\emptyset$, at least one of $u_3, u_5$ must be in $V(C)$. Also, since $C$ is $2$-connected, if there is only one of $u_3, u_5$, say $u_3$, being in $V(C)$, then $u_3$ would be a cut vertex between vertices in $\{a,b\}\backslash\{u_3\}$ and vertices in $\{x,v,u_1,u_2\}$ of $C$. Thus $\{u_3,u_5\}\in V(C)$. Similarly, we can also prove that $u_2\in V(C)$ since $\{u_2, u_5\}$ is also a cut set. Note that $u_2\in V(P)$. By Fact~\ref{fact5}, it holds $\{a,b\}\cup\{u_3, u_5\}\cup V(A(P))=\{a,b\}\cup\{u_5,v,u_1,u_2,u_3\}\subset V(C)$. Since $|V(C)|\le 5$, we have $a,b\in \{u_5,v,u_1,u_2,u_3\}$. Since $a,b\notin\{x,v,u_1,u_2\}$, it holds $\{a,b\}=\{u_3,u_5\}$.\\
Claim that: we can always pick a cycle $C'$ in $G'+ab$ of length $4$ or $5$.
If $V(P)\cap V(C)=\emptyset$, pick $C'=C$ and the claim is done. If $V(P)\cap V(C)\neq\emptyset$, then by the previous proof, $ab=u_3u_5$, which is not a non-edge in $G$ if $G'=G_2$, a contradiction. 
If $G'=G_1$, we let $C'$ be the cycle of length $4$ or $5$ in $G_2=G_1+u_3u_5=G_1+ab$.
Hence, the claim is done, which implies that $G'$ is also a $\mathcal{C}_{\{4,5\}}$-saturated graph.
However,
$$|E(G')|\le|E(G)|-5<\lceil\frac{5}{4}n-\frac{3}{2}\rceil-5= \lceil\frac{5}{4}(n-4)-\frac{3}{2}\rceil=\lceil\frac{5}{4}|V(G')|-\frac{3}{2}\rceil\mbox{,}$$
which is a contradiction by the minimality of $G$. See Figure~\ref{proof31}.
\end{proof}

\begin{figure}[h]
\centering
\subfigure{\includegraphics[width=3.74in]{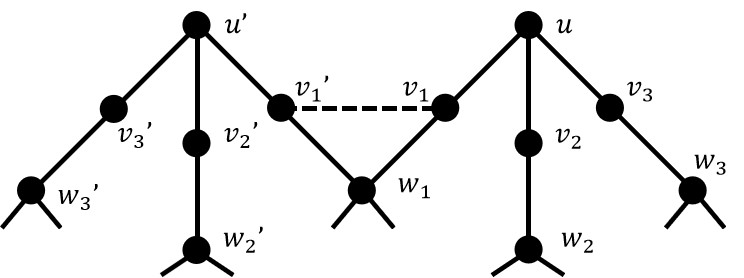}}
\caption{The proof of Lemma~\ref{only1}}\label{proof33}
\end{figure}

\begin{lem}\label{only1}
There is at most one vertex $u\in V(G)$ with $ch'(u)<0$. Moreover, $ch'(u)=-\frac{1}{4}$.
\end{lem}
\begin{proof}
Suppose to the contrary that there are two vertices $u$ and $u'$ with $ch'(u), ch'(u')<0$.
By Lemma~\ref{ch4}, $d_G(u)=d_G(u')=3$ and $N_G(u),N_G(u')\subseteq D_2^0(G)$. 
Let $N_G(u)=\{v_1, v_2, v_3\}$ and $N_G(u')=\{v_1', v_2', v_3'\}$.
Also, let $N_G(v_i)=\{u, w_i\}$ and $N_{G}(v_i')=\{u', w_i'\}$ for $i\in[1,3]$.
By Lemma~\ref{3222}, for $i\in[1,3]$, $w_i$ has a neighbor of degree at least $3$ (in $N_{G}(w_1) \cap N_{G}(w_2)\cap N_{G}(w_3)$), so $u'\neq w_i$. This means $(\{u\}\cup N_{G}(u))\cap(\{u'\}\cup N_{G}(u'))=\emptyset$.\\
Since $uu'\notin E(G)$ is a non-edge, by Fact~\ref{fact2}, there is a path $P$ of length $3$ or $4$ from $u$ to $u'$.
Since $N_G(u)=\{v_1, v_2, v_3\}$, one of $v_1, v_2, v_3$ must be in $V(P)$. Suppose $v_1\in V(P)$ without loss of generality. By Fact~\ref{fact5}, $uv_1w_1=A(v_1)\subset P$. Similarly, we can suppose $w_1'v_1'u'=A(v_1')\subset P$. 
If $w_1\neq w_1'$, then $|V(A(v_1)\cup A(v_1'))|=6>5\ge |V(P)|$, a contradiction. So $w_1=w_1'$.\\
Since $v_1v_1'\notin E(G)$ is a non-edge, by Fact~\ref{fact2}, there is a path $P'$ of length $3$ or $4$ from $v_1$ to $v_1'$.
Since $N_G(v_1)=\{u, w_1\}$, one of $v_1u$ and $v_1w_1$ is in $E(P')$. Similarly, one of $v_1'u'$ and $v_1'w_1$ is in $E(P)$.
If $v_1w_1, v_1'w_1 \in E(P')$, then $P'=v_1w_1v_1'$ is of length $2$, a contradiction. Hence, one of $v_1u$ and $v_1'u'$ is in $E(P')$. Without loss of generality, suppose $v_1u\in E(P')$.
Since $N_G(u)\backslash\{v_1\}=\{v_2, v_3\}$, one of $v_2, v_3$ must be in $V(P')$. Suppose $v_2\in V(P')$ without loss of generality. By Fact~\ref{fact5}, $uv_2w_2=A(v_2)\subset P'$. So $v_1uv_2w_2\subset P'$.\\
Now if $v_1'u'\in E(P')$, similarly, we can suppose $w_2'v_2'u'v_1'\subset P'$. However, in this case, $P'$ contains at least $6$ edges, which is a contradiction. Therefore, $w_1v_1'\in E(P')$.
Note that $w_1\neq w_2$ and $w_1w_2\notin E(G)$ by Lemma~\ref{3222}, There are at least two edges other than $v_1u, uv_2, v_2w_2, w_1v_1'$ in $E(P')$, which means $P'$ is a path of length at least $6$, a contradiction. The main part of the proof is done. See Figure~\ref{proof33}.\\
For the ``moreover" part, note that for the only vertex $u$ with $ch'(u)<0$, it is easy to see by Lemma~\ref{ch4} that $r_G(u)=3$ and $s_G^+(u)=s_G^-(u)=t_G(u)=0$. One can check by Lemma~\ref{chother} that $ch'(u)=-\frac{1}{4}$. 
\end{proof}

By the rules (R1) to (R4), $\sum_{v\in V(G)}ch(v)=\sum_{v\in V(G)}ch'(v)$.
Thus by Lemma~\ref{only1}, $\sum_{v\in V(G)}ch(v)=\sum_{v\in V(G)}ch'(v)\ge-\frac{1}{4}$. The proof of Theorem~\ref{45} is done.

\section{More on $\mathcal{C}_I$-saturated graphs}
In this paper, we considered $\mathcal{C}_I$-saturated graphs and showed that $\sat(n, \mathcal{C}_{\{4,5\}})=\lceil\frac{5n}{4}-\frac{3}{2}\rceil$.
As put in the introduction, for common positive integer $n$, $\sat(n, C_4)\neq \sat(n, \mathcal{C}_{\{4,5\}})$ (Ollmann~\cite{Oll72}, Tuza~\cite{Tuz89}, Fisher et al.~\cite{FFL97}), $\sat(n, C_5)\neq \sat(n, \mathcal{C}_{\{4,5\}})$ (Chen~\cite{Che09,Che11}) while $\sat(n, \mathcal{C}_{[4,+\infty)})=\sat(n, \mathcal{C}_{\{4,5\}})=\lceil\frac{5n}{4}-\frac{3}{2}\rceil$ (Ferrara et al.~\cite{Subdivision12}). We give the following conjecture.
\begin{conj}
For any integer $r\ge 5$, there exists a number $n(r)$, such that for any integer $n\ge n(r)$, $\sat(n, \mathcal{C}_{[4,r]})=\lceil\frac{5n}{4}-\frac{3}{2}\rceil$.
\end{conj}
Similarly, since $\sat(n, C_5)=\sat(n, \mathcal{C}_{[5,+\infty)})=\lceil\frac{10(n-1)}7\rceil$ for large enough $n$ (Chen~\cite{Che09,Che11} and Ferrara et al.~\cite{Subdivision12}). We put a similar conjecture as follows.
\begin{conj}
For any integer $r\ge 5$, there exists a number $n(r)$, such that for any integer $n\ge n(r)$, $\sat(n, \mathcal{C}_{[5,r]})=\lceil\frac{10(n-1)}7\rceil$.
\end{conj}
More generally, we give the following conjecture.
\begin{conj}
There exist functions $r(s)$ on $s\in [4,+\infty)$ and $n(s,r)$ for $s\in [4,+\infty)$ and $r\in [r(s),+\infty)$,
such that for any integers $s\ge4$, $r\ge r(s)$ and $n\ge n(r,s)$,
$\sat(n, \mathcal{C}_{[s,r]})=\sat(n, \mathcal{C}_{[s,+\infty)})$.
\end{conj}

\begin{figure}[h]
\centering
\subfigure{\includegraphics[width=1.9in]{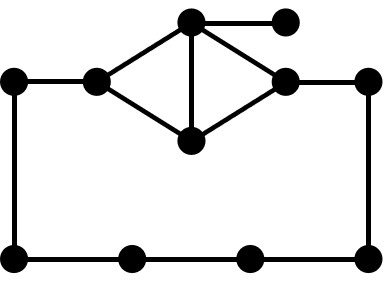}}
\caption{An example for $J_{s,t}^{+r}$ when $s=8, t=2, r=1$}\label{conj}
\end{figure}

For any integer set $I$ and integers $a$ and $b$, let $aI+b=\{ai+b: i\in I\}$. For example, $2[2,+\infty)=2\mathbb{Z}_++2=\{4,6,8,\dots\}$. 
By Observation~\ref{4cons}, one may hope that $\sat(n, \mathcal{C}_{I})=\lceil\frac{5n}{4}-\frac{3}{2}\rceil$ holds for any integer set $I$ with $\{3,4\}\cap I=\{4\}$ and $\{5,6,7\}\cap I\neq\emptyset$ when $n$ is large enough. However, this is not true. The following fact shows that $\sat(n, \mathcal{C}_{2\mathbb{Z}_++2})=n$ for any integer $n\ge 3$.
\begin{fact}
(i) If $G$ is a $\mathcal{C}_{2\mathbb{Z}_++2}$-saturated graph on $n\ge 3$ vertices, then $G$ is connected and $G$ is not a tree. In particular, $|E(G)|\ge n$ for $n\ge 3$.\\
(ii) For odd $n\ge 3$, $C_n$ is a $\mathcal{C}_{2\mathbb{Z}_++2}$-saturated graph on $n$ vertices. For even $n\ge 3$, $C_{n-1}^+$ is a $\mathcal{C}_{2\mathbb{Z}_++2}$-saturated graph on $n$ vertices, where $C_{n-1}^+$ is a graph obtained by adding a pendant edge to one vertex of $C_{n-1}$.
\end{fact}
\begin{cor}
For any integer $n\ge 3$, $\sat(n, \mathcal{C}_{2\mathbb{Z}_++2})=n$.
\end{cor}
However, we think the following conjecture may be true.
\begin{conj}
For integer $n\ge 1$, $\sat(n, \mathcal{C}_{3\mathbb{Z}_++1})=\lceil\frac{5n}{4}-\frac{3}{2}\rceil$.
\end{conj}
One may be interested in $\mathcal{C}_{a\mathbb{Z}_++b}$-saturated graphs. By the following construction, it seems to be reasonable to consider $b=2$ firstly.
For integers $s\ge 2$ and $t\ge 1$, let $J_{s,t}$ be a graph obtained by joining all the vertices of a copy of $K_t$ with the two endpoints of a copy of $P_s$. To be more specific, for integers $s\ge 2$ and $t\ge 1$, let $X=\bigcup_{i=1}^s\{x_i\}$ and $Y=\bigcup_{j=1}^t\{y_j\}$, then $V(J_{s,t})=X\cup Y$ and
$$E(J_{s,t})=E(u_1u_2\dots u_s)\cup E(K[\{u_1,u_s\},Y])\cup E(K[Y])\mbox{.}$$
For $r\in[0,t]$, let $J_{s,t}^{+r}$ be the graph obtained by adding a pendant edge to each vertex in $\{y_1,y_2,\dots,y_r\}$. 
See Figure~\ref{conj}.
Note that $C_n=J_{n-1,1}$ and $C_{n-1}^+=J_{n-2,1}^{+1}$. The following fact gives a contruction of $\mathcal{C}_{a\mathbb{Z}_++2}$-saturated graphs when $a\ge 2$.

\begin{fact}
For integers $a\ge 2$ and $n\ge a+1$, let $k\in [1,+\infty)$ and $r\in[0,a-1]$ be the unique pair with $n-1=ak+r$.
The graph $J_{a(k-1)+2,a-1}^{+r}$ is a $\mathcal{C}_{a\mathbb{Z}_++2}$-saturated graph on $n$ vertices with $n+\binom{a}{2}-1$ edges.
In particular, $\sat(n, \mathcal{C}_{a\mathbb{Z}_++2})\le n+\binom{a}{2}-1$.
\end{fact}
\begin{conj}
For integers $a\ge 2$ and $n\ge a+1$, $\sat(n, \mathcal{C}_{a\mathbb{Z}_++2})= n+\binom{a}{2}-1$.
\end{conj}

\end{document}